\documentclass[conference]{IEEEtran} 
\usepackage{fullpage,amsmath,amssymb,amsthm,citesort}
\usepackage{epsfig,graphicx,subfigure,wrapfig,psfrag}		

\def \reals {\mathbb{R}}

\def \RIPcond {\delta}			
\def \flow {\phi}				

\newcommand{\boxdim}[1]{D_{#1}}							

\newcommand{\trajectory}[1]{\widetilde{g}\left( #1 \right)}		
\newcommand{\rank}[1]{\operatorname{rank}\left( #1 \right)}		

\def \attractor {\mathcal{A}}	
\def \tattractor {\trajectory{\attractor}}	

\newtheorem{definition}{Definition}[section]
\newtheorem{lemma}{Lemma}[section]
\newtheorem{thm}{Theorem}[section]

\begin{document}

\title{A First Analysis of the Stability of Takens' Embedding}
\author{Han Lun Yap, Armin Eftekhari, Michael B. Wakin, and Christopher J. Rozell}
\maketitle

\section{Introduction}

Data in the form of time series pervades many areas of science and engineering. 
Traditionally, one might impose a suitable probability model for the data, yielding tools such as the ARMA model and the Kalman filter~\cite{book_brockwell2002timeseries}.  
More recently, modeling the time-series data as a one-dimensional measurement of the possibly high-dimensional states of a nonlinear dynamical system has also produced many positive results~\cite{kantz2004nonlinear,Asefa2005}. 

\begin{figure}
	\begin{center}
\includegraphics[width=80mm]{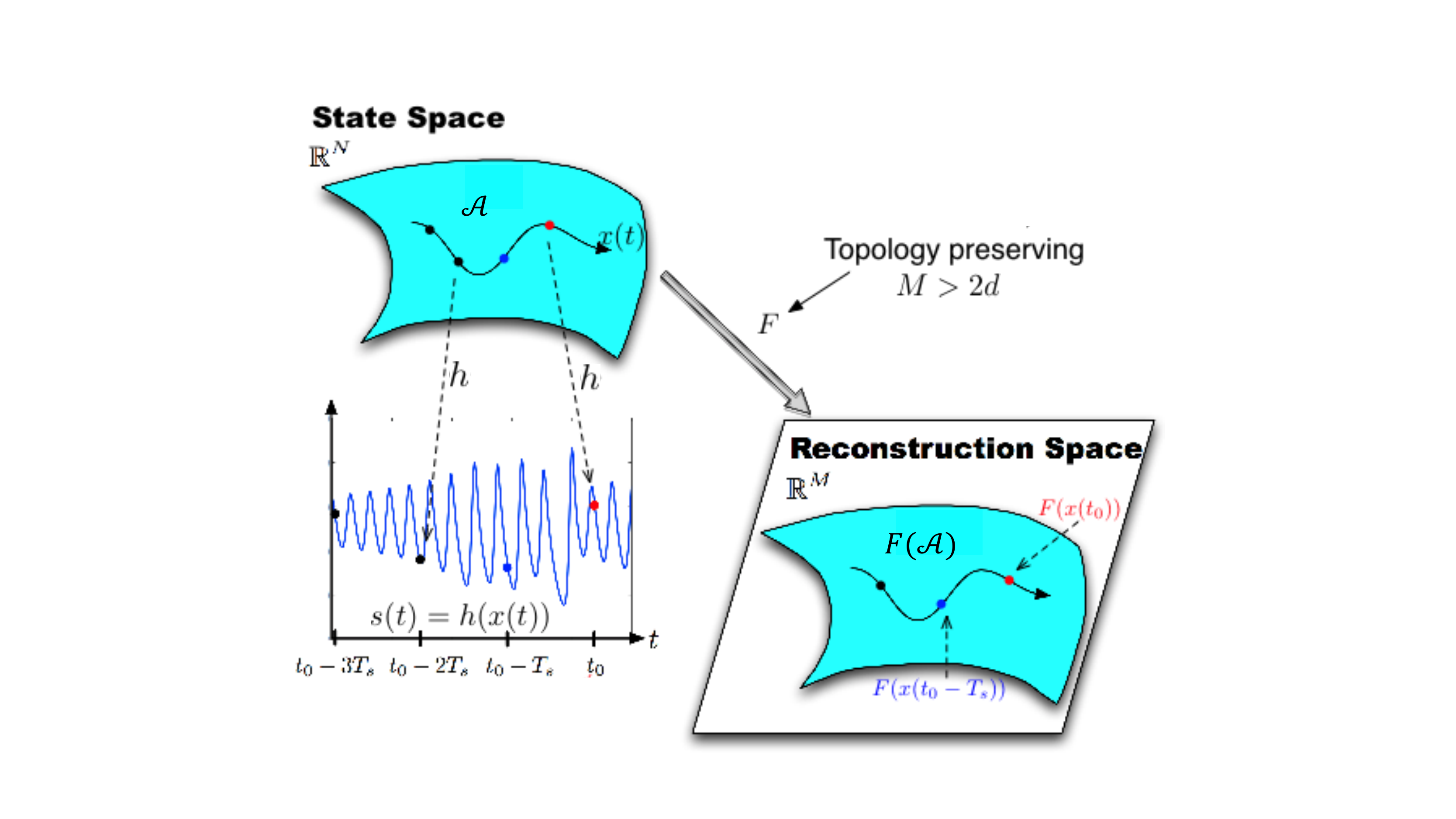}
	\end{center}
	\vspace{-5mm}
	\caption{\small Takens' embedding of a state-space attractor $\attractor$ into reconstruction space from time series measurements.}
	\vspace{-5mm}
	\label{fig:Takens}
\end{figure}

The theoretical foundation of time-series analysis based on the dynamical system model hinges on a seminal paper by Takens~\cite{Takens1981}, whose ideas were later expanded by Sauer et al.~\cite{Sauer1991}. 
The surprising result of the paper asserts that when the states of a dynamical system are confined to a low-dimensional attractor, complete information about the hidden states can be preserved in the time-series output (see Figure~\ref{fig:Takens}). 
Indeed, many systems of interest do have this type of structure~\cite{strogatz1994nonlinear}, and a variety of algorithms for tasks such as time-series prediction 
exploit Takens' result~\cite{kantz2004nonlinear}. 

To be precise, suppose we have a dynamical system whose system states $x(t)$ lie on a submanifold $\attractor \subset \reals^N$. 
Given a sampling interval $T_s$, one could define the discretized dynamics of the system through the flow function $\flow_{T_s}:\attractor \rightarrow \attractor$ defined as $x(t + T_s) = \flow_{T_s}(x(t))$. 
The time series $s(t):= h(x(t))$ is modeled as a one-dimensional observation of the system dynamics through the observation function $h$. 
Usually the sampling time $T_s$ is fixed and we will drop the subscript $T_s$ from $\flow_{T_s}$ to ease notations. 
Takens~\cite{Takens1981} defined the \emph{delay coordinate map} $F: \reals^N \rightarrow \reals^M$ as a mapping of the state vector $x(t) \in \reals^N$ of a dynamical system to a point $F(x(t)) \in \reals^M$ in the \emph{reconstruction space} by taking $M$ uniformly spaced samples of the time series $s(t)$ and concatenating them into a single vector, i.e.,
	$F(x(t)) = [s(t), \; s(t-T_s), \;\cdots, \; s(t-(M-1)T_s)]^T$.	
Takens' result states that for \emph{almost every} smooth observation function $h(\cdot)$, the delay coordinate map is an \emph{embedding} (i.e., a one-to-one immersion) of the state space attractor with dimension $d$ when $M>2d$.  
As seen in Figure~\ref{fig:Takens}, despite the state being hidden from direct observation, the topology of the attractor that characterizes the dynamical system can be preserved in the time-series data when it is arranged into a delay coordinate map. 

In the absence of measurement or system noise, Takens' result indicates that the delay coordinate vectors $F(x(t))$ are a proxy for the hidden system states $x(t)$. 
However, in the presence of noise, a one-to-one mapping as guaranteed by Takens' result may not be sufficient to guarantee the robustness of any processing  performed in the reconstruction space.  
The underlying problem is that while Takens' theorem guarantees the preservation of the attractor's \emph{topology}, it does not guarantee that the \emph{geometry} of the attractor is also preserved.  
To be precise, Takens' result guarantees that two points on the attractor do not map to the same point in the reconstruction space, but there are no guarantees that close points on the attractor remain close under this mapping (or far points remain far).  
Consequently, relatively small imperfections could have arbitrarily large, unwanted effects when the delay coordinate map is used in applications. 

In the signal processing community, recent work has highlighted the importance of well-conditioned measurement operators to ensure the geometry of a low-dimensional signal family is preserved.  
%
The notion of a \emph{stable embedding} by a measurement operator $F$ introduced in compressed sensing (CS)~\cite{Candes2006} has been extended to manifold signal families~\cite{baraniuk2009random,Clarkson2008a,Yap2012,Eftekhari2014}. 
The results in~\cite{baraniuk2009random} guarantee that distances between points on a $\boxdim{}$-dimensional submanifold are approximately preserved when $F$ is a random, linear orthoprojector and $M=O(d\log(N))$.\footnote{The notation $O(\cdot)$ simply means ``on the order of''. This result is extended in~\cite{Yap2012,Eftekhari2014} to include other random stable embedding operators. Moreover, the results in~\cite{Clarkson2008a,Eftekhari2014} removed the logarithmic dependence on $M$.} 
Besides signal recovery~\cite{Shah2011}, the stable embedding of a signal family also ensures that data processing and inference algorithms can perform with approximately the same guarantees in the measurement space as in the ambient space~\cite{Davenport2010a}. 

The conditions for Takens' Embedding Theorem to hold ignore the effects of noise. 
As such, time-series analysis in practice requires a careful empirical determination of the sampling time $T_s$ and number of delays $M$, typically via looking at the first null of mutual information or autocorrelation between two consecutive samples of the time series~\cite{fraser1986independent,kantz2004nonlinear}. 
This usually results in using a number of delay coordinates $M$ larger than the minimum prescribed by Takens' theorem~\cite{Casdagli1991a}. 

In this paper, we use tools and ideas in CS to provide a first theoretical justification for the choice of $M$ in noisy conditions. 
In fact, we show that under certain conditions on the dynamical system, measurement function $h$, number of delays $M$ and sampling time $T_s$, the delay-coordinate map can be a stable embedding of the dynamical system's attractor. 

\section{Main Results}

For the rest of this section, consider a dynamical system given by the flow $\flow = \flow_{T_s}$ with sampling time $T_s$ and whose states are lying on a Riemannian submanifold $\attractor$ of $\reals^N$. 
We shall also assume that $\flow$ is a smooth function. 
%
Before presenting the stable embedding results proper, we will need to introduce an important characterization of manifolds that will be useful for describing certain local and global properties of the system manifold $\attractor$. 
\begin{definition}
Let $\attractor$ be a compact Riemannian submanifold of $\reals^N$. The \emph{reach} $\tau$ of $\attractor$ is defined as the largest positive real number having the following property: The open normal bundle about $\attractor$ of radius $a$ is embedded in $\reals^N$ for all $a<\tau$.
\end{definition}
Reach is an equivalent measure to the condition number used in~\cite{Niyogi2009,Yap2012,baraniuk2009random}. 
See the above cites and~\cite{Eftekhari2014,federer1959curvature} for more on reach and its implications. 

On top of the properties of the system manifold, our results also depend on the types of measurement functions used to form the time series data. 
First, consider a set of basis functions $h_p : \reals^{N} \rightarrow \reals$ for $p = 1, \cdots, P$ that will span a $P$-dimensional subspace of all functions from $\reals^N$ to $\reals$. 
An element chosen from this class of measurement functions can be denoted by $h_\alpha$ for some $\alpha \in \reals^P$ defined as $h_\alpha := \sum_{p=1}^P \alpha_p h_p$. 
Next, we define the delay coordinate map with measurement function $h_\alpha$ as
	$F_{\alpha}(x) := \sum_{p=1}^{P}\alpha_{p}F_{p}(x) =\left[
	h_{\alpha}(x), \cdots, h_{\alpha}\left(\flow^{-M+1}\left(x\right)\right) \right]^T$,
where $F_p$ is defined as the delay-coordinate vector formed with the $p$-th basis function
	$F_{p}(x):= \left[h_{p}(x), \cdots, h_{p}\left(\flow^{-M+1}\left(x\right)\right) \right]^T$. 
In this paper, we are considering the class of linear functions. 
Thus, we take $P = N$ (where $N$ is the ambient dimension) and $h_p(\cdot) := \langle \cdot, e_p \rangle$ where $e_p$ is the $p$-th canonical vector. 

Additionally, we need a metric to describe the dynamical system in question. 
For any state $x \in \attractor$, we define the $M \times N$ matrix $G_x$ as
$G_{x} := \left( F_1(x) \;|\; \cdots \;|\; F_N(x) \right),$ 
and define the trajectory vector $\trajectory{x} \in \reals^{MN}$ starting at point $x \in \attractor$ as
	$\trajectory{x} := \left[ x,\flow^{-1}\left(x\right),\cdots,\flow^{-M+1}\left(x\right) \right]^T$. 
Observe that $G_{x}$ is basically the matrix version of $\trajectory{x}$ and that $F_\alpha(x) = \sum_{i = 1}^{N} \alpha_i F_i(x) = G_x \alpha$. 
We also define the trajectory manifold $\tattractor \subset \reals^{MN}$ as the set of all trajectories in $\attractor$, i.e., 
	$\tattractor:=\left\{ \trajectory{x}\,:\, x \in \attractor \right\}$. 
Notice that since the first $N$ entries of $\trajectory{x}$ are distinct for all $\trajectory{x} \in \tattractor$ and $\flow$ is a smooth function, $\tattractor$ is diffeomorphic to $\attractor$. 
Then for any $G \in \reals^{M \times N}$, define $r(G) := \frac{\|G\|_F^2}{\|G\|_2^2}$ which we shall call the \emph{soft-rank} of the matrix $G$. 
In the stable embedding results that follow, we shall see that the soft-rank of the matrices $G_x - G_y$ for all pairs $x,y \in \attractor$ capture how well information about the dynamical system lying on $\attractor$ is captured by the basis measurement functions $\{h_p\}$. 
The soft-rank is a refined metric to replace the rank of a matrix. 

Under the assumptions that there are no periodic orbits of period less than $M$ on $\attractor$ and $M \le N$, it is easy to see that $\rank{G_x - G_y} = M$ for all pairs of $x,y \in \attractor$. 
Having full rank only says that the matrix $G_x - G_y$ has $M$ non-zero singular values, but it does not say anything about the distribution of these singular values.
On the other hand, the soft-rank is the ratio of the $\ell_2$ and $\ell_\infty$ norms of the vector of singular values, therefore $r(G_x - G_y) = cM$ for some $c \in (\frac{1}{M},\; 1]$ if $G_x - G_y$ is full-rank. 
The value of $c$ tells us about the energy distribution of the singular values, with $c$ tending to 1 when all the singular values have almost the same values and tending to $\frac{1}{M}$ when only one of the singular value is significant while the rest are close to zero. 
We remark that via linear algebra, we have $c \approx 1$ whenever the rows of $G_x - G_y$ have approximately equal length and are approximately orthogonal. 

With all the required notations now defined, we finally present our stable embedding result for linear measurement functions. 
\begin{thm}
	\label{thm:manifold}
	Suppose we have a dynamical system described by a flow $\flow$ with sampling time $T_s$ and whose states lie on a $\boxdim{\attractor}$-dimensional submanifold $\attractor$ of $\reals^N$. 
	Also, suppose that the trajectory manifold $\tattractor$ has volume $V_{\tattractor}$ and reach $\tau_{\widetilde{g}(\attractor)}$. 
	Let $\alpha \in \reals^N$ be an i.i.d.\ Rademacher $\pm 1$ sequence, and let $\RIPcond$ be  a predetermined stable embedding conditioning. 
	Define the infimum soft-rank $r(\attractor):=\inf_{x \neq y \in \attractor} r(G_x - G_y)$ and assume that 
$
V_{\widetilde{g}(\mathcal{A})}\ge O\left(\tau_{\widetilde{g}(\mathcal{A})}^{D_{\mathcal{A}}}\right). 
$	
If
	\begin{eqnarray}
		r(\mathcal{A})\ge& 
		O\left(\RIPcond^{-2} D_{\mathcal{A}}\log\left( \frac{r(\mathcal{A})^{1/2} V_{\widetilde{g}(\mathcal{A})}^{1/D_{\mathcal{A}}}}{\tau_{\widetilde{g}(\attractor)}}   \right)\right), 
		\label{eq:soft_rank_cond}
	\end{eqnarray}
	then except with a probability of at most
	$
	\left(\frac{r(\mathcal{A})^{D_{\attractor}} V_{\widetilde{g}(\mathcal{A})}}{\tau_{\widetilde{g}(\attractor)}^{D_{\attractor}}}\right)^{-C}
	$
	with $C$ an absolute constant, the delay-coordinate map with $M$ delays $F_\alpha$ is a stable embedding of $\tattractor$, i.e., for all $x,y \in \attractor$,
	\begin{eqnarray}
		\label{eq:SE_statement1}
		(1 - \RIPcond) \le \frac{\left\|F_\alpha(x) - F_\alpha(y)\right\|_2^2}{\left\|\trajectory{x} - \trajectory{y}\right\|_2^2} \le (1 + \RIPcond).
	\end{eqnarray}
\end{thm}

Due to space limitations, the proof of this theorem will be deferred to a future publication. 
We note that the assumption $V_{\widetilde{g}(\mathcal{A})}\ge O\left(\tau_{\widetilde{g}(\mathcal{A})}^{D_{\mathcal{A}}}\right)$ is to maintain a failure probability less than 1. 
	

Our result describes how the geometry and the dynamics on the manifold coupled with the number of measurements $M$, as seen through the soft-rank of $G_x - G_y$ in~\eqref{eq:soft_rank_cond}, dictate the stable embedding conditioning of the delay-coordinate map. 
Thus, this result is different from the typical Takens' embedding result in two distinct ways. 
First, in the typical Takens' embedding result, only the topology of the manifold --- basically the fact that the system states lie on a manifold and the dimension of the manifold --- is utilized. 
Second, the typical result does not say what happens when we have more measurements beyond the minimum required for an embedding. 

The atypical condition~\eqref{eq:soft_rank_cond} demands further elaboration. 
Suppose the number of measurements $M$ is less than the ambient space dimension $N$ and there are no periodic orbits of period less than $M$. 
First, notice that if the \emph{infimum soft-rank} $r(\attractor)$ scales with the number of measurements $M$, then we will have a typical manifold stable embedding condition whereby a manifold is stably embedded whenever $M$ scales with the dimension of the manifold and logarithmically with other characterizations of the manifold (which is exactly what appears on the right hand side of~\eqref{eq:soft_rank_cond}), and the conditioning $\RIPcond$ of the embedding improves with increasing $M$. 
We remark that the improvement of $\RIPcond$ with $M$ matches what was observed by Casdagli et al.~\cite{Casdagli1991a}, namely that with more measurements, distortion by the delay-coordinate map is decreased. 
Now, for every $x, y \in \attractor$, the soft-rank $r(G_x - G_y)$ indeed depends on the number of measurements $M$, but it also depends on the sampling time $T_s$ coupled with various geometric properties of the manifold and the dynamical system. 
As discussed previously, the soft-rank can indeed be equal or close to $M$, but only when the rows of the matrix $G_x - G_y$ (i.e., the chords $x-y, \flow^{-1}(x) - \flow^{-1}(y), \cdots$) have almost the same length and are almost orthogonal. 
Moreover, these requirements on $G_x - G_y$ have to hold for all $x, y \in \attractor$. 

As expected, the length-preserving and orthogonality requirements described above can be very stringent.  
For example, for a dynamical system whose inverse flow $\flow^{-1}$ has a large maximal Lyapunov exponent, the requirement that the chords retain almost the same length will be easily violated for any decent sampling time $T_s$ and for $x,y$ close to one another.\footnote{
When the inverse flow $\flow^{-1}$ has Lyapunov exponent $\lambda$, we have $\|\flow^{-m}(x) - \flow^{-m}(y)\|_2 \approx e^{m\lambda}\|x - y\|_2$ for $x,y$ close enough. Therefore, the chords going down the matrix $G_x - G_y$ grow exponentially in length, severely violating the length-preserving requirement. 
}  
This is the effect of irrelevance as described in, e.g.,~\cite{Casdagli1991a}. 
Suppose $T_s$ is small and $x \approx y$. 
Then, the $M$ consecutive chords forming the $M$ rows of $G_x - G_y$ can more or less be approximated by vectors on a tangent plane of $\attractor$ (say the tangent plane of $x$). 
Since $\dim(\attractor) = \boxdim{\attractor}$ which could potentially be much less than $M$, this means that the matrix $G_x - G_y$ will have only $\boxdim{\attractor}$ significant singular values and thus $r(G_x - G_y) = cM$ with $c \ll 1$. 
This is the effect of redundancy. 
We will need many more measurements $M$ to ``break out'' of the tangent plane to achieve $r(G_x - G_y) = cM$ with a decent $c$ (i.e., $c \approx 1$). 
It is useful to add that if the manifold $\attractor$ indeed only lie on a low-dimensional subspace of $\reals^N$, then the infimum soft-rank will be \emph{upper bounded} by the dimension of this subspace and thus will \emph{not} scale with $M$. 
This remark is also true when $M$ exceeds the ambient dimension $N$, in which case the infimum soft-rank will plateau at $N$. 
Notice that this imposes the limitation that $M \le N$ for this result to be useful which thus restricts the class of dynamical systems that we can consider. 
Notably, this result is not adequate for dynamical systems residing on low-dimensional ambient spaces ($N$ small), e.g., the Lorenz system. 
Nonetheless, we intend to show in a future article that we can break this ambient space dimension ``barrier'' by considering nonlinear measurement functions.



Another important observation is that instead of embedding state-space vectors $x \in \attractor$, we are in fact embedding trajectory vectors $\trajectory{x} \in \tattractor \subset \reals^{MN}$ as seen in~\eqref{eq:SE_statement1}. 
The variables on the right-hand side of~\eqref{eq:soft_rank_cond} are also based on the geometry of $\tattractor$, e.g., the terms $V_{\tattractor}$ and $\tau_{\tattractor}$. 
Maintaining distance between trajectories of the dynamical system in the reconstruction space may be advantageous for some applications. 
We can translate this into a stable embedding of the ambient space attractor $\attractor$ since there there is a diffeomorphism (hence isomorphism) between the ambient and the trajectory manifold. 
Understandably, this translation comes with a degradation of the stable embedding conditioning $\RIPcond$ previously obtained. 
While a detailed mathematical treatment of this translation will be deferred to a future article, it suffices to note here that the requirements for the minimum degradation of the stable embedding conditioning caused by the translation usually run contrary to requirements for the best infimum soft-rank in~\eqref{eq:soft_rank_cond}.
This is also somewhat similar to the tradeoff between redundancy and irrelevancy in time series analysis. 
Redundancy happens when the sampling time and/or the Lyapunov exponent (of the inverse flow) is small. 
In this case, the trajectory and ambient vectors are scalar multiples of one another, which means that there is minimum degradation of the stable embedding conditioning going from trajectory to ambient space.  
However, the soft-rank of $G_x - G_y$ suffers because now the rows of the matrix $G_x - G_y$ are very similar to one another which violates the orthogonality condition (as described before). 
However when irrelevance occurs (meaning that either the sampling time or the Lyapunov exponent of the inverse flow becomes large), the consecutive chords in $G_x - G_y$ rapidly decorrelate. 
When this occurs, the soft-rank of $G_x - G_y$ may improve due to the diversity but the conditioning is degraded when we pass from trajectory space to ambient space. 

Lastly, we remark that our results mirror the notion that ``almost every measurement function provides an embedding'' appearing in a typical Takens' embedding statement (see \cite{Sauer1991}). 
This notion is represented by the random vector of coefficients $\alpha \in \reals^N$. 
Because $F_\alpha$ is a stable embedding with high probability on the coefficients $\alpha$, this means that most measurements functions in the space of functions defined by $\{h_\alpha = \alpha \;|\; \alpha \in \{-1, +1\}^N \}$ can result in a stable embedding by delay coordinate maps. 
In fact, the class of measurement functions considered can be vastly expanded by replacing the Rademacher sequence $\alpha$ with a Gaussian sequence instead. 
However due to the lack of space, this result will be discussed in a future article. 

\section{Experiment}
We shall demonstrate our theoretical results with a simple example. 
Consider the following discrete-time dynamical system with system states at time step $n$, $x_n \in \reals^N$, defined through
	$x_{n+1} = \Phi x_n$ 
where the (linear) flow $\Phi \in \reals^{N \times N}$ is a shift matrix, i.e.,
\begin{eqnarray*}
	\Phi = \left( 
	\begin{array}{ccccc}
		0 & 1 & 0 & \cdots & 0 \\
		\vdots &   & \ddots & \ddots & 0 \\
		0 & & & 0 & 1 \\
		1 & 0 & \cdots & \cdots & 0
	\end{array}
	\right).
\end{eqnarray*}
Also, suppose the initial condition of the system is given by $x_0 = [1, 0 \cdots, 0]^T$. 
This system can be thought of as the dynamics of a single point object translating down the entries of a $N$-dimensional vector. 
Thus, the system states $x_n$ lie on a one-dimensional manifold $\attractor$ that is parameterized by the object's location on the vector.  
We suppose that we only get to observe a one-dimensional time series derived by taking a linear projection of the system states, i.e., we get to observe the time series
	$s_n = \alpha^T x_n$ 
where $\alpha \in \reals^N$ is an i.i.d.\ Rademacher sequence. 
The goal is to study how well the delay coordinate vectors formed using the time series $s_n$ provide a stable embedding of the trajectory vectors $\trajectory{x_n}$ of the dynamical system. 

The discussion following Theorem~\ref{thm:manifold} suggests that the stable embedding conditioning is dependent on how well the infimum soft-rank $r(\attractor)$ scales with the number of measurements $M$. 
For $x \neq y \in \attractor$, $G_x - G_y$ is a circulant matrix since $\Phi$ is a shift matrix. 
Moreover, the first row of of $G_x - G_y$ is a vector of zeros except for a `$1$' and a `$-1$'. 
With these observations, the following lemma shows that the infimum soft-rank of this matrix can be calculated analytically. 
\begin{lemma}
	For the discrete-time dynamical system described in the text, we have $r(\attractor) \ge \frac{M}{2}$. 
\end{lemma}
The proof of this lemma follows from the fact that $(G_x - G_y) (G_x - G_y)^T$ is also a circulant matrix and the formulas for the eigenvalues of circulant matrices are well-known. 
This lemma tells us that for this system, we get a linear scaling of the infimum soft-rank with $M$. 
Thus Theorem~\ref{thm:manifold} will imply that the stable embedding conditioning $\RIPcond$ will scale with $M^{-1/2}$. 



\section{Conclusion}

This paper is a first step towards understanding whether Takens' embedding can be stablilized through concatenating more time series measurements. 
Many important issues are not yet addressed in this work. 
These include the study of nonlinear measurement functions, the translation of the results from trajectory manifold to state-space manifold, and the extension to fractal attractors instead of just Riemannian manifold attractors. 
We intend to address these issues in a future work. 

\bibliographystyle{IEEEtran}
\bibliography{ref}

\end{document}